\documentclass[11pt,twoside,a4paper]{article}
\usepackage{amsmath,amssymb,amsthm}

\usepackage{graphicx,psfrag}

\newtheorem{Thm}{Theorem}[section]
\newtheorem{Lem}[Thm]{Lemma}
\newtheorem{Pro}[Thm]{Proposition}

\theoremstyle{definition}

\theoremstyle{remark}
\newtheorem{Rem}[Thm]{Remark}
\newtheorem{Def}[Thm]{Definition}

\newcommand{\R}{\mathbb{R}}

\newcommand{\N}{\mathbb{N}}
\newcommand{\Hy}{\mathbb{H}}

\newcommand{\al}{\alpha}

\newcommand{\ga}{\gamma}

\newcommand{\de}{\delta}

\newcommand{\om}{\omega}

\newcommand{\la}{\lambda}

\renewcommand{\phi}{\varphi}

\renewcommand{\th}{\theta}

\newcommand{\ca}{\operatorname{ca}}

\renewcommand{\d}{\partial}
\newcommand{\di}{\d_{\infty}}

\newcommand{\sub}{\subset}

\newcommand{\ov}{\overline}

\newcommand{\wh}{\widehat}

\begin{document}

\title{Hyperbolicity, ${\rm CAT(-1)}$-spaces and the
Ptolemy Inequality }
\author{Thomas Foertsch  
\ \& Viktor Schroeder\footnote{Supported by Swiss National
Science Foundation}}

\date{}
\maketitle

\begin{abstract} 
Using a four points inequality for the boundary of
${\rm CAT(-1)}$-spaces we study the relation between Gromov
hyperbolic spaces and ${\rm CAT(-1)}$-spaces.

\end{abstract}

\section{Introduction}

From various hyperbolic cone constructions it is known that every bounded, complete 
metric space can appear as the visual boundary of a Gromov hyperbolic space. Here visual 
boundary means the boundary of a Gromov hyperbolic space endowed with a visual metric.
In order to study the relation of (rough geodesic) Gromov hyperbolic spaces and $\operatorname{CAT}(-1)$-spaces
in terms of asymptotic methods, it is of major importance to understand which metric spaces can
appear as visual boundaries of $\operatorname{CAT}(-1)$-spaces. Surprisingly enough, due to our knowledge 
there does not appear any necessary condition for this in the literature. One of the main purposes of 
this paper is to provide a first such condition, namely a
four point relation, which we will call the {\em Ptolemy Inequality}:

\begin{Thm}\label{thm:ptolemy}
Let $Y$ be the boundary of a ${\rm CAT(-1)}$-space endowed 
with a Bourdon or a Hamenst\"adt metric
$| \ |$. Let
$y_1,y_2,y_3,y_4 \in Y$, then

$$|y_1y_3||y_2y_4| \le |y_1y_2||y_3y_4| + |y_2y_3||y_4y_1|.$$

Equality holds if and only if the convex hull of
the four points is isometric to an ideal quadrilateral in the 
hyperbolic plane $\Hy^2$ such that the geodesics $y_1y_3$ and
$y_2y_4$ are the diagonals.
\end{Thm}

Note that the formulation of this four point inequality
is M\"obius invariant. Thus, if the metric 
$|\ |$ is replaced by a M\"obius equivalent metric, the inequality is invariant.
Therefore it holds for all Bourdon metrics and also for all
Hamenst\"adt metrics on $Y$ (for a discussion of these metrics
compare Section \ref{subsec:catminusone}).
As a consequence the inequality is well adapted to the geometry of
the boundary of a ${\rm CAT(-1)}$-space. 
It is a classical theorem attributed to Ptolemy (85-165), that
if $y_1,\ldots ,y_4$ are points in this order on a circle
in the Euclidean plane, then we have equality in this formula.
The classical Ptolemy theorem is equivalent to the {\em if direction} of the
equality discussion in Theorem \ref{thm:ptolemy}.
We will obtain the inequality from a detailed study
of the proof of a result of Bourdon \cite{B2}.

We use 
Theorem \ref{thm:ptolemy} to study the relation
between Gromov hyperbolic spaces and
${\rm CAT(-1)}$-spaces.
Clearly every 
${\rm CAT(-1)}$-space is Gromov hyperbolic.
Since Gromov hyperbolicity is not a local curvature
condition, the opposite is not true in general.
Given a Gromov hyperbolic space $X$ 
one can ask the following question:
Is $X$ rough
isometric to some
${\rm CAT(-1)}$-space $W$? 
Here a map $f:X \to W$ between metric spaces is
called a {\em rough isometric embedding}, if there exists a
constant $R \ge 0$ such that 
$$|xx'| - R \le |f(x)f(x')| \le |xx'| + R.$$
If in addition the image
$f(X) \sub W$ is $R$-dense, then $f$ is called a 
{\em rough isometry}.

We look for answers to this question for
Gromov hyperbolic spaces which are
{\em visual} (see Section \ref{sec:hypspace}).
The visual condition can be viewed as a quasiisometry-invariant
version of
the condition of extendable geodesics.
For simplicity of the exposition the reader may think that 
$X$ is a geodesic space with a basepoint
$o \in X$ such that for every point $x \in X$ there exists
$y \in \di X$ such that 
$x$ lies on a geodesic $oy$.

Coming back to our question we remark that
Gromov hyperbolicity is invariant under arbitrary scaling
of the metric while the 
${\rm CAT(-1)}$-condition is only invariant under scaling
with factors $\la \le 1$.
Therefore the formulation of the problem is not yet good enough
for our purposes.

If $(X,d)$ is a metric space, we can consider
the whole family of scaled metric spaces $(X, \la\, d)$
with $\la >0$ and look for a distinguished
normalization. We use the asymptotic upper
curvature bound 
$K_u(X)$ defined in \cite{BF} as a normalization
(see Section \ref{sec:hypspace}). 

This normalization is only possible if
$K_u(X)$ is finite. In the case that
$K_u(X)= - \infty$, the space $X$ looks very much like a tree.
We call a visual Gromov hyperbolic space
{\em treelike} if $K_u(X) = -\infty$.
This definition is justified by the result
in \cite{BF}, that a visual Gromov hyperbolic space
with $K_u(X)= -\infty$ is rough isometric to a tree provided that
in addition $\di X$ is doubling.
Since trees are ${\rm CAT(-1)}$,
it is not a substantial restriction to consider only
nontreelike spaces.

\smallskip

Thus in the sequel we will consider only
nontreelike visual Gromov hyperbolic spaces.
Let $X$ be such a space, then we can normalize
$X$ such that $K_u(X)= -1$.
We call the normalized metric on $X$ the critical metric and
use the symbol 
$d_0$ for it.

\smallskip

{\bf Main Question}:
Let $(X,d_0)$ be a nontreelike visual Gromov hyperbolic space 
endowed with its
critical metric $d_0$. 
Does there exist a ${\rm CAT(-1)}$ space $W$, such that
$(X,d_0)$ is rough isometric to $W$?

\smallskip

In this setting one can reformulate the Bonk-Schramm embedding result (compare \cite{BS}).
It says that under the doubling condition on $\di X$ the answer to this
question is {\em almost yes} in the sense that one only needs an
arbitrarily small scaling of the critical metric to obtain
the desired rough isometry. More precisely the Bonk-Schramm result 
(which relies on the  Assouad
embedding theorem) implies the
following:

\begin{Thm} \label{thm:bonkschramm}
Let $(X,d_0)$ be a nontreelike visual Gromov hyperbolic space endowed
with its critical metric.
Assume in addition that the boundary $\di X$ is doubling.
Then for every positive  $\la < 1$ there exists
a rough isometry
of 
$(X,\la\, d_0)$ to a 
${\rm CAT(-1)}$-space $W$.
\end{Thm}

\begin{Rem}
Actually in \cite{BS} it was proven more explicitly that 
for every positive  $\la < 1$ there 
exists a number $N$,
such that 
$(X,\la\, d_0)$ is rough isometric to a 
convex subset $W$ of the
standard hyperbolic space
$\Hy^N$.
\end{Rem}

We want to remark that
it follows from the definition of the critical metric,
that $(X,\la d_0)$ cannot be rough isometric to a 
${\rm CAT(-1)}$ space for any $\la >1$
(compare Remark \ref{rem:critexp}). 

A related embedding result can be
obtained by combining
results of Lang-Schlichenmaier (see \cite{LS})
and
Alexander-Bishop (see \cite{AB}).

\begin{Thm} \label{thm:lsab}
Let $(X,d_0)$ be as above and assume now that the boundary
$\di X$ has finite Assouad Nagata dimension. Then there exists
some $\la < 1$ such that
$(X,\la\, d_0)$ is rough isometric to a
${\rm CAT(-1)}$-space.
\end{Thm}

Details for a proof of this theorem will be given elsewhere. Just note that 
by a theorem of Lang-Schlichenmeier every metric space of finite
Assouad Nagata dimension admits a snowflake embedding into a product of 
a finite number of metric trees. 
This product certainly is a ${\rm CAT (1)}$-space.
Now Alexander-Bishop construct ${\rm CAT (-1)}$-spaces as certain metric warped 
products with fibers that are ${\rm CAT (1)}$-spaces. 
In order to establish the
validity of Theorem \ref{thm:lsab}, it only remains to verify that the 
fiber's ${\rm CAT (1)}$-metric actually yields a visual metric on the boundary
at infinity of such a ${\rm CAT (-1)}$ warped product. 
The embedding statement 
of Theorem \ref{thm:lsab} then just follows exactly as the one of Theorem
\ref{thm:bonkschramm}. \\

Although the proof of Theorem \ref{thm:lsab} needs some
constant $\la$ bounded away from $1$ (more precisely: the proof of
Theorem 1.3 in \cite{LS} needs this constant), it is unknown if 
the result is true for any positive
$\la < 1$ (similar as in the Bonk-Schramm-Assouad result).

One main result of our paper is the existence of
an example of a Gromov hyperbolic space
$(X,d_0)$ such that the optimal $\la$ for which 
$(X,\la\, d_0)$ is rough isometric to a
${\rm CAT(-1)}$-space is bounded away from $1$.
This implies in particular that the main question
as stated above has a negative answer.

\begin{Thm} \label{theo-main}
There exists a visual Gromov hyperbolic space 
$(X,d_0)$ with the following property. If
$\frac{1}{2}<\la$, then there does not exist a ${\rm CAT(-1)}$-space
$W$ which is rough isometric to $(X,\la\,d_0)$. 
However,
$(X,\frac{1}{2}\,d_0)$ is rough isometric to a
${\rm CAT(-1)}$-space.
\end{Thm}

Our result allows now to reformulate the question more
quantitatively and to introduce a new invariant
$\la_0$ for visual nontreelike Gromov hyperbolic spaces:

{\bf Questions:} Let $(X,d_0)$ be a visual Gromov hyperbolic space with
its critical metric. Does there exist some $0< \la \le 1$ such
that $(X,\la\, d_0)$ is rough isometric to a
${\rm CAT(-1)}$-space? If such $\la$ exists, what is the supremum $\la_0$ of
these $\la$? Is $(X,\la_0\,d_0)$ rough isometric to some ${\rm CAT(-1)}$-space?

\begin{Rem}
Note that the set 
$\la$, such that $(X,\la\,d_0)$ is rough isometric to some
${\rm CAT(-1)}$-space is  either empty or an interval of the form
$(0,a)$ or $(0,a]$ with $a\le 1$.
\end{Rem}


\noindent We give an outline of the paper.
Let $X$ be a Gromov hyperbolic space.
We denote 
by $Z = \di X$
the boundary at infinity of $X$.
Given a basepoint
$o \in X$, the expression
$e^{-(.|.)_o}$ defines a quasi-metric
on $Z$, here $(.|.)_o$ denotes the Gromov product.

If $o$ and $o'$ are different basepoints then the
quasi-metrics
$e^{-(.|.)_o}$ and
$e^{-(.|.)_{o'}}$ are 
bi-Lipschitz. 
Thus the bi-Lipschitz class
$[\rho]$ of the quasi-metric
$\rho =e^{-(.|.)_o}$ is well defined
and does not depend on
the basepoint.



If we scale the metric on $X$ by a factor
$\la$,
then the Gromov product
$(.|.)_o$ is transformed into
$\la (.|.)_o$ and the corresponding 
quasi-metric on $Z$ is taken to the power $\la$.
Thus it is reasonable to consider the whole family
$\rho^{\la}$ of quasi-metrics and not only the
particular quasi-metric $\rho=e^{-(.|.)_o}$.

Given a general quasi-metric space
$(Z,\rho)$ we can associate to $\rho$ a
critical exponent
$s_0 \in (0,\infty]$ 
(see Section \ref{subsec:quasimetric}).
If $s_0 \neq \infty$, we say that
$\rho^{s_0}$ is the critical quasi-metric on $Z$.
In the case that $X$ is a visual Gromov hyperbolic space,
consider $Z = \di X$ endowed with the quasi-metric
$\rho = e^{-(.|.)_o}$. Then there is a relation 
of the critical exponent $s_0$ of $\rho$ and the
asymptotic upper curvature bound $K_u(X)$ defined in
\cite{BF}. Indeed it holds
$K_u(X) = - s_0^2$. 

If $X$ is nontreelike (i.e. $s_0 \ne \infty$), then
one can scale the metric on $X$ in a unique way, such that
$e^{-(.|.)_o}$ (where now the Gromov product is taken with respect
to the scaled metric) is in the critical class. This
corresponds to the scaling 
$K_u(X)=-1$. In this way we find a distinguished metric
$d_0$ on $X$.

We are interested in the question, if one can embed $X$ 
rough isometrically into some
${\rm CAT(-1)}$-space $W$.
The existing embedding theorems work in the following way.
First find an embedding of the boundary
$\di X= Z$ into the boundary of some
${\rm CAT(-1)}$-space, i.e.
a map
$f: Z \to Y$, where $Y$ is the boundary of
some $W$. Then one extends this embedding
to an embedding
$F:X \to W$.
The idea of the extension is
easily 
explained in the case
that $X$ is a geodesic Gromov hyperbolic space
with extendable geodesics.
Given a basepoint $o\in X$ and an arbitrary point
$x \in X$, there exists a point
$z \in Z$ and a geodesic
$oz$, such that $x \in oz$.
The extension $F$ is now defined as follows.
Choose some basepoint
$o' \in W$. Now define
$F(x)$ to be the point on the geodesic
$o'f(z)$ such that
$|ox|=|o'F(x)|$.

Bonk and Schramm proved that
$F$ is a rough isometric embedding
if and only if
$f$ is a bi-Lipschitz map.
In this case $F$ is a rough isometry onto
the convex hull of $F(X)\subset Y$, which turns
out to be ${\rm CAT(-1)}$ itself.

Using this extension construction,
the embedding problem can be reduced to
an embedding problem
$f:Z \to Y$,
where $Z$ is some complete bounded quasi-metric
space, and $Y$ is the boundary of a
${\rm CAT(-1)}$-space
(endowed with a Bourdon metric).

To discuss this embedding problem, we recall here the definition of
a snowflake map.
A map
$f:Z \to Y$ between quasi metric spaces is called a
{\it $q$-snowflake map}, if
there exists $c > 1$ such that for all $z,z'\in Z$
$$\frac{1}{c} |zz'|^q \le |f(z)f(z')| \le c|zz'|^q.$$

This means that the quasi-metric
$\rho^q$ embeds bi-Lipschitz into the metric
space $Y$. In the case
that $\rho$ is critical, we conclude in particular that
$q\leq 1$ and $q=1$ can only occur, if
the critical quasi-metric
$\rho$ is actually bi-Lipschitz to a metric.

Thus we have the following:
Let $X$ be a nontreelike visual Gromov hyperbolic space, then
$(X,\la\,d_0)$ can be rough isometrically
embedded into a ${\rm CAT(-1)}$-space,
if there exists
a $\la$-snowflake map
from $\di X$ to the boundary $Y$ of a 
${\rm CAT(-1)}$-space.

To obtain our example we
denote with $Z$ the unit ball in $\ell^1$,
i.e. a point $z \in Z$ is a sequence
$(z_1,z_2,\ldots )$ with
$\sum |z_i|\le 1$. We prove

\begin{Thm} \label{thm:snowflaketheorem}
For
$q>\frac{1}{2}$
there does not exist a $q$-snowflake 
embedding
$f:Z \to Y$ 
where $Y$ is a space
satisfying the
Ptolemy inequality.
\end{Thm}

However we show there exists a $\frac{1}{2}$-snowflake map
of $Z$ into some Hilbert space, which
is the boundary of the infinite dimensional
hyperbolic space (see Section \ref{sec:example}). \\

Finally note that Theorem \ref{thm:snowflaketheorem} is sharp with respect to our methods of proof
in the following sense: What we actually prove is that $Z$ does not admit a $q$-snowflake
embedding into a metric space satisfying the Ptolemy inequality for $q>\frac{1}{2}$. However,
in Section \ref{sec-ptolemy} we prove 
\begin{Pro} \label{prop-sqrt-ptolemy}
Let $(X,d)$ be an arbitrary metric space. Then $(X,d^{1/2})$ satisfies the Ptolemy inequality.
\end{Pro}
Thus, if we want to obtain a non embedding theorem similar in spirit to Theorem \ref{thm:snowflaketheorem}
with snowflake parameters $q\le \frac{1}{2}$, then we need other necessary conditions for a metric space to
appear as Bourdon or Hamenst\"adt metrics on the boundary of $\operatorname{CAT} (-1)$-spaces. \\
We would like to emphazise that the most natural candidate allowing such a non embedding result seems to 
be the unit ball in $l^\infty$. \\

It is a pleasure to thank Sergei Buyalo, Mario Bonk, Alexander Lytchak
and Urs Lang for many discussions about hyperbolic spaces.

\section{Preliminaries}

\subsection{Quasi-metrics and metrics} \label{subsec:quasimetric}

A {\em quasi-metric space}
is a set
$Z$
with a function
$\rho:Z\times Z\to [0,\infty)$
which satisfies the conditions:

\noindent (1) $\rho(z,z')\ge 0$
for every
$z$, $z'\in Z$
and
$\rho(z,z')=0$
if and only if
$z=z'$;

\noindent (2)
$\rho(z,z')=\rho(z',z)$
for every
$z$, $z'\in Z$;

\noindent (3)
$\rho(z,z'')\le K\max\{\rho(z,z'),\rho(z',z'')\}$
for every
$z$, $z'$, $z''\in Z$
and some fixed
$K\ge 1$.

Let 
$(Z,\rho)$ be a quasi-metric space.
By $[\rho]$ we denote the bi-Lipschitz class of $\rho$,
i.e. for a map
$\rho' : Z\times Z \to [0,\infty)$ we have
$\rho' \in [\rho]$ if and only if there exists $c \ge 1$ such that
for all $z,z' \in Z$
$$\frac{1}{c} \ \rho(z,z') \le \rho'(z,z') \le c\ \rho(z,z').$$

We are interested in obtaining a metric
on $Z$.
Since the only problem is the triangle inequality, the following 
approach is very natural.
Define a map
$d:Z\times Z \to [0,\infty]$,
$d(z,z')=\inf\sum_i\rho(z_i,z_{i+1}),$
where the infimum is taken over all sequences
$z=z_0,\dots,z_{n+1}=z'$
in
$Z$.
By definition $d$ satisfies the
triangle inequality. 
We call this approach to the triangle
inequality the 
{\it chain approach}.

For a quasi-metric $\rho$ we denote
with
$d =\ca(\rho)$ the pseudometric which we obtain
when applying the chain approach to $\rho$.

The problem with the chain approach is that
$d(z,z')$ could be $0$ for different
points $z, z'$ and axiom
(1) is not longer satisfied 
for $(Z,d)$.

Frink \cite{Fr} realized that the chain approach works
for $2$-quasi-metric spaces. 

\begin{Pro}\label{pro:2quasimetric}  Let
$\rho$
be a
$2$-quasi-metric
on a set
$Z$ and let for 
$z,z' \in Z$,
$d(z,z')=\inf\sum_i\rho(z_i,z_{i+1}),$
where the infimum is taken over all sequences
$z=z_0,\dots,z_{n+1}=z'$
in
$Z$.
Then $d$ is a metric on $Z$ with
$\frac{1}{4}d \le \rho\le d$.

\end{Pro}

If $(Z,\rho)$ is a quasi-metric space,
then $\rho^s$ is a 2-quasi-metric if
$s>0$ is sufficiently small.

\begin{Def}
A quasi-metric space $(Z,\rho)$ is called {\em LM-space}
(Lipschitz metrizable),
if $\ca(\rho) \in [\rho]$.

\end{Def}

\noindent Hence a quasi-metric space is LM if and only if the following
two conditions hold:

\noindent (1) the chain approach gives a metric.

\noindent (2) the metric from the chain approach is bi-Lipschitz to
$\rho$.

\noindent Clearly the LM property is a bi-Lipschitz
invariant.

One easily proves the following:
If $\rho$ is LM, then $\rho^s$ also is LM for
every $0<s\le 1$.

Note that $\rho^s$ is a $2$-quasi-metric for $s$ small enough.
Thus to every quasi-metric space $(Z,\rho)$ which is not
bi-Lipschitz to an ultrametric one can associate in a unique way
a {\em critical exponent} $s_0 \in (0,\infty]$ with the property:
$\rho^s$ is LM for all $s< s_0$ and $\rho^s$ is not LM for all $s > s_0$.

\begin{Rem}\label{rem:critexp} 
It follows from the definition of the critical exponent, that
for every $s > s_0$ there cannot exist a bi-Lipschitz map
$(Z,\rho^s) \to Y$, where $Y$ is a metric space.
\end{Rem}

We shortly discuss the situation $s_0=\infty$.
It follows from \cite{BF}:

\begin{Thm} \label{thm:critex}
Let $(Z,\rho)$ be a doubling quasi-metric space.
Then $s_0 = \infty$ 
if and only if $\rho$ is bi-Lipschitz to an ultrametric .
\end{Thm}

Without the doubling assumption Theorem \ref{thm:critex} fails in general. This 
follows from an example due to Leonid Kovalev: $Consider$ the set of integers $\mathbb{N}$ 
endowed with the metric $d$, where $d(m,n):=\log (1+|m-n|)$. For this metric the critical
exponent $s_0$ is $\infty$, but $(\mathbb{N},d)$ is not bi-Lipshitz equivalent to an ultrametric.
Nevertheless, we call (ad hoc) a quasi-metric $\rho$ {\em ultrametriclike} if $s_0=\infty$.


\subsection{Gromov hyperbolic spaces} \label{sec:hypspace}

Let $X$ be a metric space.
For $o,x,x' \in X$ let
$$(x|x')_o := \frac{1}{2} ( |ox| + |ox'| - |xx'|). $$
The space $X$ is called {\em $\de$-hyperbolic} if for 
$o,x,x',x'' \in X$
\begin{equation} \label{eq:de-eq1}
((x|x')_o,(x|x'')_o,(x'|x'')_o)
\ \ \mbox{is a} 
\ \de\mbox{-triple}
\end{equation}
in the sense that the two smallest of the three numbers differ by at most
$\de$.

$X$ is called {\em hyperbolic}, if it is
$\de$-hyperbolic for some
$\de \ge 0$.
The relation (\ref{eq:de-eq1}) is called the
{\em $\de$-inequality} with respect to the point $o\in X$.

If $X$ satisfies the $\de$-inequality for one
individual
basepoint $o \in X$, then it satisfies
the $2\de$-inequality for any other basepoint
$o' \in X$ (see for example \cite{G}). Thus, to check
hyperbolicity, one has to check this inequality
only for one basepoint.

Let
$X$
be a hyperbolic space and
$o\in X$
be a base point. A sequence 
$\{x_i\}$ of points $x_i\in X$
{\em converges to infinity,}
if
$$\lim_{i,j\to\infty}(x_i|x_j)_o=\infty.$$
Two sequences
$\{x_i\}$, $\{x_i'\}$
that converge to infinity are {\em equivalent}
if
$$\lim_{i\to\infty}(x_i|x_i')_o=\infty.$$
Using the
$\de$-inequality,
one easily sees that this defines an equivalence relation
for sequences in
$X$
converging to infinity. The {\em boundary at infinity}
\index{boundary at infinity}
$\di X$
of
$X$
is defined as the set of equivalence classes
of sequences converging to infinity.

For points
$y$, $y'\in\di X$
we define their Gromov product by
$$(y|y')_o=\inf\liminf_{i\to\infty}(x_i|x_i')_o,$$
where the infimum is taken over all sequences
$\{x_i\}\in y$, $\{x_i'\}\in y'$.
Note that
$(y|y')_o$
takes values in
$[0,\infty]$ and
that
$(y|y')_o=\infty$
if and only if
$y=y'$.
In a similar way we define for $\xi \in \di X$,
$x\in X$
$$(y|x)_o=\inf\liminf_{i\to\infty}(x_i|x)_o.$$

If $X$ is 
$\de$-hyperbolic and
if $y,y',y'' \in \di X $ ,
then
$((y|y' )_o,(y |y'')_o,(y' |y''))$ is
a $\de$-triple.
This implies that the expression
$\rho(y,y') = e^{-(y|y' )_o}$ defines a
$K$-quasi-metric on $\di X$ where
$K = e^{\de}$.

A Gromov hyperbolic space is called
{\em visual}, if there exists a point $o\in X$ and a constant
$D \ge 0$ such that for every $x \in X$ there exists $y \in \di X$ with
$|ox| - (x|y)_o \le D$.

Roughly speaking, in a visual Gromov hyperbolic space the position
of a point $x$ is (up to a universal constant), given
by some point $y \in \di X$ and the distance $|ox|$ from the
basepoint. It turns out that for these spaces almost all information is encoded
in the properties of
$\di X$.

On the other hand, if some bounded metric space $Y$ is given,
then it is possible to construct a Gromov hyperbolic space
$X$ such that $\di X$ as a set coincides with $Y$ and
the quasi-metric
$e^{-(.|.)_o}$ is bi-Lipschitz to the given metric
on $Y$ (see for example \cite{BS}).

Also the following holds.
Let $X$ be a visual Gromov hyperbolic space.
Then $\di X$, endowed with the quasi-metric
$\rho = e^{-(.|.)_o}$, is bi-Lipschitz to an ultrametric if and only if
$X$ is rough isometric to a tree.

In \cite{BF} the notion
$K_u(X)$ of an upper asymptotic curvature bound is introduced.
In the case of visual Gromov hyperbolic spaces this notion is 
strongly related to the critical exponent of the quasi-metric
$e^{-(.|.)_o}$ on $\di X$.
The following relation holds (see Theorem 1.5 in \cite{BF}):
$K_u(X)= - s_0^2$.

A visual Gromov hyperbolic space is called {\em treelike}, if
$K_u(X) = -\infty$.

This definition is motivated by the following result (see \cite{BF}):

\begin{Thm}
Let $X$ be a visual Gromov hyperbolic space. Assume in addition
that $\di X$ is doubling.
Then $K_u(X)=-\infty$ if and only if $X$ is rough isometric to a tree.
\end{Thm}



\subsection{${\rm CAT(-1)}$ spaces}\label{subsec:catminusone}

Let now $X$ be a ${\rm CAT(-1)}$ space,
i.e. $X$ is a complete geodesic metric space,
such that triangles are thinner than comparison triangles
in the hyperbolic plane
$\Hy^2$. In particular $X$ is also Gromov hyperbolic. Let
$Y= \di X$.
Given $x \in X$ and $w \in X \cup \di X$ there exists a unique
geodesic segment $xw$ from $x$ to $w$.
If $y_1,y_2 \in \di X$ are different points, there is also
a unique geodesic line $y_1y_2$ joining these points.

Given a point $o \in X$ and points
$y_1,y_2 \in \di X$
we denote by
$\angle_o(y_1,y_2)$ the local angle at
$x$, i.e. the angle between the initial directions of
the geodesics from $oy_1$ and $oy_2$.
By
$\th_o(y_1,y_2)$ we denote the asymptotic comparison angle.
I.e. let $y_i(t)$ be the point on the ray $oy_i$ with
distance $t$ to $o$. Let 
$\ov{o}, \ov{y_1(t)}, \ov{y_2(t)}$ be the comparison triangle in
$\Hy^2$, and let $\ov{\ga_t}$ be the angle of this triangle at
$\ov{o}$. Then $\th_o(y_1,y_2)=\lim_{t\to\infty}\ov{\ga_t}$.
Let
$\rho_o(y_1,y_2) = \sin(\frac{1}{2}\th_o(y_1,y_2))$ 
be the Bourdon metric (with basepoint $o$).
Indeed Bourdon proved \cite{B1} that this expression is a metric
and satisfies the triangle inequality.
One can also express $\rho_o$ in terms of the Gromov product
and obtains the formula
$\rho_o(y_1,y_2) = e^{-(y_1|y_2)_o}$, i.e. in the ${\rm CAT(-1)}$ situation
the Bourdon metric corresponds to the quasi-metric considered earlier.
For the convenience of the reader we give a proof of
this formula. The computation is in the hyperbolic plane
$\Hy^2$. The triangle $\ov{o}, \ov{y_1(t)}, \ov{y_2(t)}$ has a
limit ideal triangle given by
geodesic rays
$\ga_i:[0,\infty)\to\Hy^2$
starting from
$\ov{o}$ with angle
$\theta_o$. We then
have
$$e^{-(y_1|y_2)_o} = \lim_{t\to\infty} (e^{h_t}e^{-2t})^{1/2},$$
where
$h_t=d(\ga_1(t),\ga_2(t))$
is the distance in
$\Hy^2$.
From the hyperbolic law of cosine
$$\cosh(h_t) = \cosh^2(t) - \sinh^2(t)\cos\theta_o$$
and the trigonometric formula
$1-\cos\theta_o=2\sin^2(\theta_o/2)$,
we easily obtain
$$e^{h_t}\sim e^{2t}\sin^2(\theta_o/2)$$
as
$t\to\infty$.
Hence, the claim.

Bourdon metrics with respect to different basepoints
are M\"obius equivalent.
Thus 
given a fixed Bourdon metric $|\ |$ on
$Y$, we have for an arbitrary $o\in X$ that
$$\frac{\rho_o(y_1,y_2)\rho_o(y_3,y_4)}{\rho_o(y_1,y_3)\rho_o(y_2,y_4)}
=\frac{|y_1y_2||y_3y_4|}{|y_1y_3||y_2y_4|}.$$
The M\"obius invariance can be seen as follows.
Let $o' \in X$ be a different basepoint, then a trivial
computation shows for $x_1,x_2 \in X$ that
$$(x_1|x_2)_{o'}=|oo'|+(x_1|x_2)_o -(x_1|o')_o-(x_2|o')_o,$$
a formula which extends to points $y_1,y_2 \in \di X$. Thus
$$\rho_{o'}(y_1,y_2)=\mu \frac{\rho_o(y_1,y_2)}{\la(y_1)\la(y_2)},$$
where $\mu=e^{-|oo'|}$ and $\la(y)=e^{-(y|o')_o}$. This clearly implies
that $\rho_{o'}$ is M\"obius equivalent to $\rho_o$.

We should mention here that 
Hamenst\"adt \cite{H} introduced (even earlier) a metric on
$\di X \setminus \{\om\}$, where $X$ is a Hadamard manifold with
curvature $\le -1$ and $\om \in \di X$ is a distinguished point.
We do not describe her construction verbatim but modify her
construction such that it works also for general
${\rm CAT(-1)}$ spaces:
fix a point
$\om \in \di X$ and consider a Busemann function
$b$ for the point $\om$. Define the Gromov product with respect
to this Busemann function, i.e.
$$(x|x')_b =\frac{1}{2}(b(x)+b(x')-|xx'|),$$
which also extends to points at infinity.
The corresponding {\it Hamenst\"adt metric} $\rho_b$ is then defined
as
$e^{-(.|.)_b}$.
If $b$ is the Busemann function at $\om$ such that
$b(o)=0$ for some point $o\in X$, then one easily computes
$b(x)=|ox|-2(\om|x)_o$ and by straightforward calculation one obtains the formula
$$(x|x')_b = (x|x')_o - (\om|x)_o -(\om|x')_o,$$
which also extends to infinity. Hence
$$\rho_b(y,y') = \frac{\rho_o(y,y')}{\rho_o(y,\om) \rho_o(y',\om)}.$$
Thus the Hamenst\"adt metric can be obtained by involution
at the point $\om$ from the Bourdon metric and in particular these
metrics are M\"obius equivalent.

Note that by definition
$\rho_b$ is only a quasi-metric on
$\di X\setminus \{\om\}$, since the involution of an arbitrary metric
does not necessarily satisfy the triangle inequality.
However in our situation the triangle inequality
$$\rho_b(y,y'') \le \rho_b(y,y') + \rho_b(y',y'')$$
is equivalent to the Ptolemy inequality
$$\rho_o(y,y'')\rho_o(y',\om) \le
\rho_o(y,y')\rho_o(y'',\om) +\rho_o(y',y'')\rho_o(y,\om),$$
which is proved in the next section. Thus
$\rho_b$ is actually a metric. This observation can be considered
as the first application of the Ptolemy inequality. As a sideremark of
this observation we formulate this in larger generality:

\begin{Rem}
Let $(Z,d)$ be an arbitray metric space. For $z \in Z$ consider
the involution
$d_z:Z \setminus \{z\} \times Z \setminus \{z\} \to [0,\infty)$,
$d_z(a,b)=d(a,b)/(d(a,z)d(b,z))$. 
Then $d_z$ is a metric for all $z \in Z$ iff $d$ satisfies the
Ptolemy inequality.

\end{Rem}

Observe that the ${\rm CAT(-1)}$ condition implies that
$\angle_o(y_1,y_2)\le \th_o(y_1,y_2)$.
However the following holds:
if
$\th_o(y_1,y_2)= \pi$ then $\angle_o(y_1,y_2)=\pi$.
($\th_o(y_1,y_2)= \pi$ implies 
that $d(y_1(t),y_2(t)) =2t$ and hence
$\angle_o(y_1,y_2)=\pi$.)

Thus we conclude that
$\rho_o(y_1,y_2)\le 1$ with equality if and only if
$o$ lies on the geodesic $y_1y_2$.

We will use the following

\begin{Lem} \label{lem:sincomp}
Let $X$ be a ${\rm CAT(-1)}$ space,
$y_1,\ldots,y_4\in \di X$ be different points and $o\in X$, then
$$\sin\frac{1}{2}\th_o(y_1,y_2)\sin\frac{1}{2}\th_o(y_3,y_4)
\le \frac{|y_1y_2||y_3y_4|}{|y_1y_3||y_2y_4|}.$$
Equality holds if and only if
the geodesics from $y_1y_3$ and $y_2y_4$ intersect in the point
$o$.

\end{Lem}

\begin{proof}
We have
\begin{eqnarray}
\sin\frac{1}{2}\th_o(y_1,y_2)\sin\frac{1}{2}\th_o(y_3,y_4)
& = & \rho_o(y_1,y_2)\rho_o(y_3,y_4) \nonumber \\
& \le & \frac{\rho_o(y_1,y_2)\rho_o(y_3,y_4)}{\rho_o(y_1,y_3)\rho_o(y_2,y_4)}
= \frac{|y_1y_2||y_3y_4|}{|y_1y_3||y_2y_4|}. \nonumber
\end{eqnarray}
Equality holds if and only if
$\rho_o(y_1,y_3)=\rho_o(y_2,y_4)=1$, which is equivalent
to $o$ lying on
$y_1y_3\cap y_2y_4$.

\end{proof}

\section{The Ptolemy Inequality}
\label{sec-ptolemy}

We start this section with a proof of Proposition \ref{prop-sqrt-ptolemy}, which states that
the square root of any metric space satisfies the Ptolemy inequality. The following proof is
due to Alexander Lytchak: \\
{\bf Proof of Proposition \ref{prop-sqrt-ptolemy}:} 
Let $\{ x_1,y_1,x_2,y_2\}$ be an ordered quadrupel in $(X,d)$ and set $p_1:=d(x_1,x_2)$,
$p_2:=d(y_1,y_2)$, $q_1:=d(x_1,y_1)$, $q_2:=d(x_2,y_1)$, $q_3:=d(x_2y_2)$ and $q_4:=d(x_1,y_2)$.
Without loss of generality we may assume that
\begin{displaymath}
\begin{array}{lclclc}
p_1 & \le & q_1+q_2 & \le & q_3+q_4 & \mbox{and} \\
p_2 & \le & q_2+q_3 & \le & q_1+q_4. &
\end{array}
\end{displaymath}
Our claim follows, once we verify that $\sqrt{p_1p_2}\le \sqrt{q_1q_3} \sqrt{q_2q_4}$. We prove this
inequality by showing that for suitable $p_1'\ge p_1$, $p_2'\ge p_2$ and $q_4'\le q_4$ one obtains
$\sqrt{p_1'p_2'}\le \sqrt{q_1q_3}\sqrt{q_2q_4'}$. \\
The numbers $p_1'$, $p_2'$ and $q_4'$ are obtained as follows. First
set $p_1':=q_1+q_2$ and $p_2':=q_2+q_3$. Then choose $q_4'$ such that
\begin{eqnarray}
q_3+q_4'=q_1+q_2 & \mbox{and} & q_2+q_3\le q_1+q_4' \hspace{1cm} \mbox{or} \label{eqn-rel1}\\
q_3+q_4'\le q_1+q_2 & \mbox{and} & q_2+q_3=q_1+q_4'.\nonumber
\end{eqnarray}
Without loss of generality we may assume that the relations (\ref{eqn-rel1}) are satisfied.
It follows that $q_1-q_3=q_4-q_2$ and $-(q_1-q_3)\le q_4-q_2$ from which we deduce $q_4\ge q_2\ge 0$
and $\epsilon :=q_1-q_3\ge 0$. Now the inequality
\begin{displaymath}
2q_2q_3 \; \le \; 2\sqrt{q_3(q_3+\epsilon )q_2(q_2+\epsilon )}
\end{displaymath}
implies
\begin{displaymath}
\begin{array}{crcl}
& \sqrt{(q_2+q_3+\epsilon )(q_2+q_3)} & \le & \sqrt{q_3(q_3+\epsilon )} \; + \; \sqrt{q_2(q_2+\epsilon )} \\
& & & \\
\Longleftrightarrow & \sqrt{p_1'p_2'} & \le & \sqrt{q_1q_3} \, \sqrt{q_2q_4'},
\end{array}
\end{displaymath}
from which the claim follows.
\begin{flushright}
$\Box$
\end{flushright}

Next we prove 
\begin{Thm}\label{thm:comparison}
Let $Y$ be the boundary of a ${\rm CAT(-1)}$-space with a Bourdon or
Hamenst\"adt metric
$|.|$. Let
$y_1,y_2,y_3,y_4 \in Y$, then

$$|y_1y_3||y_2y_4| \le |y_1y_2||y_3y_4| + |y_2y_3||y_4y_1|.$$

Equality holds if and only if the convex hull of
the $y_i$ is isometric to an ideal quadrilateral in $\Hy^2$, such that
the geodesics $y_1y_3$ and $y_2y_4$ intersect.

\end{Thm}

\begin{Rem}
We note that this is a M\"obius invariant comparison statement.
I.e. if the equality is true for some metric $|.|$, then it also
true for every M\"obius equivalent metric.
Thus this kind of comparison result is suitable for the boundary of a
${\rm CAT(-1)}$-space.
\end{Rem}

{\bf Proof of Theorem \ref{thm:comparison}:}
We use an idea from the proof of
Lemma 3.1 in Bourdon's paper \cite{B2}.

Consider the geodesic
$y_2y_4$ in our ${\rm CAT(-1)}$-space.
By continuity there exists a point
$x\in y_2y_4$ with
$\theta_x(y_1,y_2)=\theta_x(y_3,y_4)$.
Denote this angle by $\beta$.
Thus, by
Lemma \ref{lem:sincomp}  we obtain the following estimate
$$\sin^2\frac{1}{2}\beta \le
\frac{|y_1y_2||y_3y_4|}{|y_1y_3||y_2y_4|}.$$
Let 
$\ga=\theta_x(y_2,y_3)$ and
$\delta =\theta_x(y_4,y_1)$,
then the same argument shows
$$\sin\frac{1}{2}\gamma \sin\frac{1}{2}\delta
\le
\frac{|y_2y_3||y_4y_1|}{|y_1y_3||y_2y_4|}.$$
Since $x \in y_2y_4$ we see
$\angle_x(y_2,y_3) + \angle_x(y_3,y_4) \ge \pi$ and
$\angle_x(y_4,y_1) + \angle_x(y_1,y_2) \ge \pi$.
Since $\angle_x \le \theta_x$ we obtain
$\beta + \ga \ge \pi$ and
$\beta + \delta \ge \pi$. Consequently 
$\sin \frac{1}{2} \ga \ge \sin \frac{1}{2} (\pi - \beta)
= \cos \frac{1}{2} \beta$ and also
$\sin \frac{1}{2} \delta \ge \cos \frac{1}{2} \beta$.
It follows
$$\sin\frac{1}{2}\gamma \sin\frac{1}{2}\delta \ge 
\cos^2 \frac{1}{2}\beta,$$ 
and hence
$$\cos^2 \frac{1}{2}\beta \le \frac{|y_2y_3||y_4y_1|}{|y_1y_3||y_2y_4|}.$$
Since
$\cos^2 \frac{1}{2}\beta +\sin^2\frac{1}{2}\beta =1$,
we obtain the Ptolemy inequality.

If we have equality, then
we have 
actually
equality in all estimates.
This implies now the rigidity statement by the following arguments.
By Lemma \ref{lem:sincomp} the diagonals
$y_1y_3$ and $y_2y_4$ intersect at the point
$x$. Furthermore we have equality for the angles
$\angle_x(y_1,y_2)=\theta_x(y_1,y_2)$,
$\angle_x(y_2,y_3)=\theta_x(y_2,y_3)$,
$\angle_x(y_3,y_4)=\theta_x(y_3,y_4)$,
$\angle_x(y_4,y_1)=\theta_x(y_4,y_1)$.
This implies by standard rigidity results
that the ideal triangles
$xy_iy_{i+1}$ are isometric to triangles in
$\Hy^2$. Since, moreover, the angles at $x$ add up to
$2\pi$ we finally see that the span of
$y_1,\ldots ,y_4$ is isometric to an ideal quadrilateral in
the hyperbolic plane. \hfill $\Box$

{\bf Examples:} 
At the end of this section we
give two examples. We show that any four point metric space
which satisfies the Ptolemy inequality can be M\"obius embedded into
the boundary of some ${\rm CAT(-1)}$-space.
Then we give an example
of a metric space (of six points) satisfying the Ptolemy inequality but
which cannot be M\"obius embedded into
the boundary of a ${\rm CAT(-1)}$-space.

Let us first consider a four point metric
space
$W=\{w_1,\ldots ,w_4\}$ with metric $d$ which satisfies
the Ptolemy inequality.
Let
$a_1=d(w_1,w_2)$, $a_2=d(w_3,w_4)$,
$b_1=d(w_2,w_3)$, $b_2=d(w_4,w_1)$,
$c_1=d(w_1,w_3)$, $c_2=d(w_2,w_4)$.
We define a new metric $d'$ on $W$ by
setting
$a'_1=a'_2=a'=\sqrt{a_1a_2}$,
$b'_1=b'_2=b'=\sqrt{b_1b_2}$,
$c'_1=c'_2=c'=\sqrt{c_1c_2}$.

Any nonsingular triangle in $(W,d')$ has the distances
$a',b',c'$ and since the numbers
$a'^2,b'^2,c'^2$ satisfy the triangle inequality by
the Ptolemy inequality, also
$a',b',c'$ satisfies the triangle inequality.
We can (after renumbering the points) assume that 
$c'\ge \max \{a',b'\}$.
Finally we scale the metric and let
$c=1$, $a=a'/c'$, $b=b'/c'$.

This new metric is M\"obius equivalent to $d$ and thus we can
start with this metric. The Ptolemy inequality now says
$a^2+b^2 \ge 1$.
Define the angles $\al$ and $\beta$ such that
$\sin \frac{1}{2} \al = a$ and
$\sin \frac{1}{2} \beta = b$. Then 
$a^2+b^2\ge 1$ implies
$\al + \beta \ge \pi$.
Let $\Delta_{\al}$ be the ideal triangle $oyy'$ in $\Hy^2$, such
that $\angle_o(y_,y')=\al$. Now glue four triangles
$\Delta_{\al} \Delta_{\beta}\Delta_{\al}\Delta_{\beta}$ 
cyclically together to obtain an ideal quadrilateral
formed from hyperbolic pieces with cone angle
$2(\al+\beta) \ge 2\pi$. Thus the space is ${\rm CAT(-1)}$,
and by construction the Bourdon metric at the cone point
coincides with our given metric.

To obtain the six point example
we start with some general remark.
Let $X$ be a ${\rm CAT(-1)}$-space and let
$Y = \di X$. Let $|.|$ be some fixed Bourdon or Hamenst\"adt
metric on $Y$.

Assume that there are points
$y_1, \ldots ,y_4 \in Y$ such that we have
equality in the Ptolemy inequality.
Then by the equality case above the geodesics $y_1y_3$ and $y_2y_4$ intersect
in some point $x \in X$.
Let $a, b$ be the positive numbers such that
$$a^2 = \frac{\rho_x(y_1,y_2)\rho_x(y_3,y_4)}{\rho_x(y_1,y_3)\rho_x(y_2,y_4)} =
\frac{|y_1y_2||y_3y_4|}{|y_1y_3||y_2y_4|},$$
$$b^2 = \frac{\rho_x(y_2,y_3)\rho_x(y_4,y_1)}{\rho_x(y_1,y_3)\rho_x(y_2,y_4)} =
\frac{|y_2y_3||y_4y_1|}{|y_1y_3||y_2y_4|}.$$
Note that $\rho_x(y_1,y_3)=\rho_x(y_2,y_4)=1$, since
$x$ lies on the corresponding geodesics.
Since we have equality in the Ptolemy inequality,
we see that $a^2 + b^2 = 1$.
Let
$\al \in (0,\pi)$ such that
$$\sin^2 \frac{1}{2}\al = \rho_x(y_1,y_2)\rho_x(y_3,y_4)= a^2,$$
$$\cos^2 \frac{1}{2}\al = \rho_x(y_2,y_3)\rho_x(y_4,y_1)= b^2.$$
Then by the above calculations
$\al$ and $(\pi - \al)$ is the intersection angle of the
geodesics $y_1y_3$ and $y_2y_4$. 
An immediate consequence is:

\begin{Lem}
Let $y_1,\ldots ,y_4$ be points in $Y$ with
equality in the Ptolemy inequality.
Assume in addition that
$|y_1y_2||y_3y_4|=|y_2y_3||y_4y_1|$.
Then the four points span a hyperbolic quadrilateral
and the geodesics $y_1y_3$ and $y_2y_4$
intersect orthogonaly.
\end{Lem}

To construct our example
consider first the standard hyperbolic
3-space $\Hy^3$.
Let $o\in \Hy^3$ and consider three
geodesics through $o$ which are pairwise
orthogonal.
We denote by
$ e_i^{\pm} \in \di \Hy^3$, $i=1,2,3$, the endpoints
of these geodesics. Then the Bourdon metric $\rho_o$ satisfies
$\rho_o(e_i^+,e_i^-)=1$ for $i=1,2,3$ and all other distances are equal to
$\frac{1}{\sqrt{2}}$.

We now define on this 6 point space $\{e_i^{\pm}| i=1,2,3\}$ another metric
$|.|$.
Therefore choose $a,b,c \in \R$ numbers close to $1$.
Define $|e_i^+ e_i^-| = 1$ for $i=1,2,3$.

$$|e_1^+ e_2^+ |=|e_1^+ e_2^- | = \frac{a}{\sqrt{2}},
\ \ |e_1^- e_2^+ |=|e_1^- e_2^- | =  \frac{1}{a \sqrt{2}},$$

$$|e_2^+ e_3^+|=|e_2^+ e_3^- | = \frac{b}{\sqrt{2}},
\ \ |e_2^- e_3^+|=|e_2^- e_3^-| = \frac{1}{b \sqrt{2}},$$

$$|e_3^+ e_1^+|=|e_3^+ e_1^-| = \frac{c}{\sqrt{2}},
\ \ |e_3^- e_1^+|=|e_3^- e_1^-| =  \frac{1}{c \sqrt{2}}.$$

If the numbers $a,b,c$ are close to $1$, then the triangle
inequality is still satisfied.
Note that for the three four-point-configurations
$\{ e_1^{\pm}, e_2^{\pm}\}$, $\{e_1^{\pm}, e_3^{\pm}\}$,
$\{e_2^{\pm}, e_3^{\pm}\}$ we have still equality
in the corresponding Ptolemy inequality.
One easily checks that all other four-point-configurations 
still satisfy the Ptolemy inequality if
$a,b,c$ are close to $1$.

We claim that the metric space
$(\{e_i^{\pm}\},|.|)$ cannot be M\"obius embedded into the boundary 
of a ${\rm CAT(-1)}$-space.
Assume to the contrary that it embeds into $\di X$.
Then by the above Lemma the three geodesics
$e_i^+e_i^-$ intersect pairwise orthogonaly.
Let us first assume
that the three geodesics intersect in one common point
$x \in X$. Then the above Lemma together with the fact, that
each pair of geodesics spans a hyperbolic quadrilateral implies
that the metric $\rho_x$ coincides with the metric
$\rho_o$ in $\Hy^3$. But note that (for $a,b,c$ different from 1),
the metic $\rho_o$ is not M\"obius equivalent to
the metric $|.|$: e.g.
$$\frac{|e_1^+ e_2^-||e_2^+e_3^+|}{|e_1^+e_2^+||e_2^-e_3^+|}= b^2 \ne 
1 = \frac{\rho_o(e_1^+,e_2^-)\rho_o(e_2^+,e_3^+)}{\rho_o(e_1^+,e_2^+)\rho_o(e_2^-,e_3^+)}.$$

Thus the three geodesics do not intersect in a common point. Hence we obtain three
intersection points $x_{12}, x_{23},x_{31}$ of the corresponding geodesics 
which form a triangle with three right angles.
Since the space is ${\rm CAT(-1)}$, such a triangle cannot exist. \\

Finally note that every $\operatorname{CAT}(0)$-space satisfies the Ptolemy inequality.
This follows from the fact that the euclidean plane satisfies the Ptolemy inequality
(which is classical), and the fact that every four point configuration in a 
$\operatorname{CAT}(0)$-space admits a subembedding into the euclidean plane
(compare \cite{BH}, p. 164).


\section{Proof of Theorem \ref{thm:snowflaketheorem}}

The proof is inspired by results of Enflo
(compare \cite{E1} and \cite{E2}).
Let $Y =\di X$, where $X$ is a ${\rm CAT(-1)}$ space.
On $Y$ we consider a Bourdon metric.

\subsection{Cubes in spaces satisfying the Ptolemy inequality}

An $m$ cube in a metric space is a subset of
$2^m$ points which are indexed by the set
$\{0,1\}^m$. 
Thus a $2$-cube in $Y$ are four points
$y_{(0,0)}$, $y_{(0,1)}$, $y_{(1,1)}$, $y_{(1,0)}$.
On the set of indices we consider the
Hamming metric $d_H$, i.e. the distance between two
indices is the number of different entries.
The points $y_I$, $I\in \{0,1\}^m$ are called vertices.
A pair of points $y_I$, $y_J$ is called a
$d$-diagonal, if $d_H(I,J)=d$. The $1$-diagonals are also called
{\it sides}. The distance
$|y_Iy_J|$ is called the length of the diagonal.

We denote by
$S_{n,m}=\{I \in \{0,1\}^n|d_H(I,0)=m\}$.
Thus $S_{n,m}$ is the set of 
$\{ 0,1\}$-sequences of length $n$ containing exactly
$m$ entries $1$.
Note that
$d_H(I,J)$ is even for
$I,J \in S_{n,m}$.

We will first consider certain homothetic embeddings
of the cube
$\{0,1\}^m$ into $S_{n,m}$ where $n \ge 2m$.
For $1\le i < j \le n$ consider the map
$\phi_{i,j}:\{0,1\}\to \{0,1\}^n$ defined by
$\phi_{i,j}(0) = e_i$,
$\phi_{i,j}(1) = e_j$.
For a sequence
$1\le k_1<k_2<\ldots <k_{2m}\le n$ we define
$\phi_{k_1\cdots k_{2m}}:\{0,1\}^m\to S_{n,m}$ by
$\phi_{k_1\cdots k_{2m}}(i_1,\ldots ,i_m)=
\phi_{k_1k_2}(i_1)+\ldots +\phi_{k_{2m-1}k_{2m}}(i_m)$.
For example consider
$\phi_{1245}:\{0,1\}^2 \to S_{5,2}$, which maps

\noindent $(0,0)\mapsto (1,0,0,1,0)$;
$\ \ \ (0,1)\mapsto (1,0,0,0,1)$;
$\ \ \ (1,1)\mapsto (0,1,0,0,1)$;
$(1,0)\mapsto (0,1,0,1,0)$.

These maps are homotheties with factor 2, i.e.
for every multiindex
$K=k_1\cdots k_{2m}$ and for all
$I,J \in \{0,1\}^m$ we have
$d_H(\phi_K(I),\phi_K(J) = 2d_H(I,J)$.

\begin{Thm} \label{theo-cubes}
For every $m \in \N$ there exists some
$n=n_m \ge 2m$ with the following property:
Let $Y$ be 
a metric space satisfying the Ptolemy inequality
and let
$\Phi:S_{n,m} \to Y$ be a map 
into $Y$, such that 
$|\Phi(I)\Phi(J)|\le b$ if $d_H(I,J)=2$ for some $b> 0$ (recall that 2 is the
minimal nontrivial distance).
 Then there exists a multiindex
$K$ with $1\le k_1<k_2<\ldots <k_{2m}\le n$, such that
the map
$\Phi\circ \phi_K:\{0,1\}^m \to Y$ is a map of an
$m$-cube into $Y$ such that there exists a diagonal 
in $\{0,1\}^m$ the image of which under $\Phi \circ \phi_K$
has length $\le  \sqrt{m}\ b$.
\end{Thm}

\begin{proof}
The proof is by induction on $m$, where the case $m=1$
is trivial.
Let us assume that the result is true for 
the value $m-1$, with $n_{m-1}$ being the corresponding
$n$-value.
Let
$$p= 2^{m-1} {{n_{m-1} }\choose{2m-2}} + 1.$$
We define
$n=n_{m-1}+p$ and will show now that the result is true for
$n= n_m$.
Let therefore
$\Phi: S_{n,m} \to Y$ be a map as in the statement of the theorem.
For $1\le i\le p$ we define canonical embeddings
$\rho_i:S_{n_{m-1},m-1} \to  S_{n,m}$ by
$\rho_i(I)= (I,0,\ldots,0,1,0\ldots,0)$, where we put
the entry $1$ at the $i$'th additional place,
i.e. at the $(n_{m-1}+i)$'th place of the sequence.
We now apply the induction hypothesis to the
$\Phi_i= \Phi\circ \rho_i:S_{n_{m-1},m-1} \to Y$.
Thus for every $i$ there exists a multiindex
$K_i$ such that the maps
$\tau_i=\Phi \circ \rho_i \circ \phi_{K_i}$ satisfy
the requirement of the statement. 
There are only ${n_{m-1}}\choose{2m-2}$ different multi-indices.
By the choice of $p$ we see that there exists a common multiindex
$K'$ with 
$1\le k_1<k_2<\ldots <k_{2m -2}\le n_{m-1}$ such that
at least
$2^{m-1} +1$ of the maps
$\tau_i=\Phi \circ \rho_i \circ \phi_{K'}$ satisfy
the requirement of the statement.
Since the cube
$\{0,1\}^{m-1}$ has $2^{m-1}$ diagonals, there are
at least two of the $i$'s
(lets call them $i_1<i_2$), where the same diagonal 
in the image has length $\le \sqrt{m-1} \ b$. 
Let $a, \ov{a} \in \{0,1\}^{m-1}$ be the endpoints of the diagonal.
Define
$k_{2m-1}=n_{m-1}+i_1$ and
$k_{2m}=n_{m-1}+i_2$ and consider
the multiindex
$K$ which is obtained from $K'$ by
extending it through $k_{2m-1}k_{2m}$.
We show that the map
$\Phi\circ \phi_K :S_{n,m}\to Y$ contains a diagonal of
length $\le \sqrt{m}$.

Note that restricted to
$\{0,1\}^{m-1}\times\{0\}$ this map coincides with
$\tau_{i_1}$ and restricted to
$\{0,1\}^{m-1}\times\{1\}$ with
$\tau_{i_2}$.
We consider the images
$y_{(a,0)}$, $y_{(a,1)}$, $y_{(\ov{a},0)}$,
$y_{(\ov{a},1)}$ of the corresponding four points in $\{0,1\}^m$.
We have 
$|y_{(a,0)}y_{(\ov{a},0)}|,
|y_{(a,1)}y_{(\ov{a},1)}|
 \le  \sqrt{m-1}\ b$ by induction hypothesis.
Furthermore we have
$|y_{(a,0)}y_{(a,1)}|,
|y_{(\ov{a},0)}y_{(\ov{a},1)}|
\le b$ by assumption on $\Phi$.
By Theorem \ref{thm:comparison} this implies that the product
$$|y_{(a,0)}y_{(\ov{a},1)}|\ |y_{(\ov{a},0)}y_{(a,1)}|
\le (m-1)b^2 +  b^2 =  m b^2,$$
which implies that at least one of the diagonals
$y_{(a,0)}y_{(\ov{a},1)}$ or
$y_{(a,1)}y_{(\ov{a},0)}$ has length
$\le  \sqrt{m}\ b$.
\end{proof}

We now prove Theorem \ref{thm:snowflaketheorem}.

Assume that there exists a $(q,c)$-snowflake map
$\Psi:Z\to Y$, where $Z$ is the unit ball in $\ell^1$.
Let $m \in \N$ be given and let $n=n_m$ as above.
Consider the map
$s_m:S_{n,m} \to Z$, 
$I \mapsto \frac{1}{m} I$. Here we consider
$\frac{1}{m} I$ as a sequence in $\ell^1$ (with
m entries of value $\frac{1}{m}$ and all other values 
$0$). Thus $s_m$ is just a scaling map the image of which lies
on the unit sphere of $Z$.
We consider
$\Phi=\Psi \circ s_m$ and apply Theorem \ref{theo-cubes}
to it.
Note that for $I,J \in S_{m,n}$ with $d_H(I,J)=2$ we
have
$$|\Phi(I)\Phi(J)| \le c 2^q \frac{1}{m^q},$$
since
$d_Z(s_m(I),s_m(J))=2\frac{1}{m}$ and the snowflake property  
of $\Psi$. 
If $I, J$ have distance $d_H(I,J) =2m$, then
$d_Z(s_m(I),s_m(J))=2$ and hence
$$|\Phi(I)\Phi(J)| \ge \frac{2^q}{c}$$
by the snowflake property.
By Theorem \ref{theo-cubes} there exists a diagonal, i.e. 
points $I,J$ with $d_H(I,J) =2m$, such that
$$|\Phi(I)\Phi(J)|\le \sqrt{m} c 2^q\frac{1}{m^q}.$$
Thus 
$\sqrt{m} \frac{1}{m^q}\ge \frac{1}{c^2}$ which implies
(since we can choose $m$ arbitrarily large and independent from
$c$) that $q \le 1/2$.


\section{An example} \label{sec:example}

Let $Z$ be as above the unit ball in
$\ell^1$.
We construct in this section a
$\frac{1}{2}$-snowflake map
$f:Z \to Y$, where $Y$ is the boundary of some
${\rm CAT(-1)}$-space $X$.
Therefore we first give a $\frac{1}{2}$-snowflake embedding of
$Z$ into the Hilbert space $\ell^2$.

By the Assouad embedding theorem there exists
$N \in \N$ and a bi-Lipschitz embedding
$(\R,d^{1/2}) \to \R^N$, i.e.
a $\frac{1}{2}$ snowflake map
$h:\R \to \R^n$.
Thus
there is some constant $c$, such that
$$\frac{1}{c} |t-s| \le |h(t)h(s)|^2 \le c |t-s|,$$
where we consider the Euclidean metric on $\R^N$.
Now the map
$g:Z \to \ell^2$,
$(z_1,z_2,\ldots )\mapsto (h(z_1),h(z_2),\ldots)$
satisfies also
$$\frac{1}{c} |zz'| \le |g(z)g(z')|^2 \le c |zz'|,$$
for all $z,z' \in Z$.
Thus by $g$ the space $Z$ is mapped into a bounded ball
of a Hilbert space and $g$ is a $\frac{1}{2}$-snowflake map.

To map it to the boundary of a
${\rm CAT(-1)}$ space, we consider as
$X$ the infinite dimensional version of
the hyperbolic space in the unit ball model.
The boundary $Y$ is the unit sphere in a Hilbert
space.
The relation is as in the classical Euclidean situation
with
the classical stereographic
projection
$\phi:S^n\to\wh\R^n$,
$$\phi(x)=\frac{1}{1-x_{0}}(x_1,\dots,x_n)\
  \text{for}\ x=(x_0,\dots,x_{n}).$$
  
Here we consider
$S^n \sub \R^{n+1}$, where on 
$\R^{n+1}$ we have coordinates
$x=(x_0,\dots,x_{n})$.

\begin{figure}[htbp]
\centering
\psfrag{o}{$o$}
\psfrag{i}{$e_0$}
\psfrag{sn}{$S^n$}
\psfrag{whrn}{$\wh\R^n$}
\psfrag{sri}{$S_r(e_0)$}
\includegraphics[width=0.4\columnwidth]{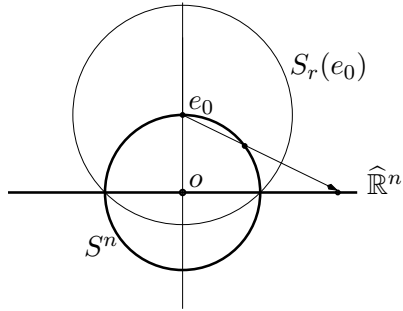}
\caption{The stereographic projection as an inversion}\label{fi:stereo}
\end{figure}

The inversion
$\wh\phi:\wh\R^{n+1}\to\wh\R^{n+1}$
of the extended
$\wh\R^{n+1}=\R^{n+1}\cup\infty$
with respect to the sphere
$S_r(e_0)\sub\R^{n+1}$,
$e_0=(1,\dots,0)$, $r=\sqrt 2$,
restricted to
the standard unit sphere
$S^n\sub\R^{n+1}$,
coincides with the stereographic projection,
$\wh\phi|S^n=\phi$.
Thus
$\phi$
as well as its inverse
$\pi:\wh\R^n\to S^n$
are M\"obius maps.

We put
$o=(0,\dots,0)\in\R^{n+1}$
and denote by
$\rho$
the standard metric on
$\R^{n+1}$, $\rho(x,y)=|x-y|$,
canonically extended to
$\wh\R^{n+1}$.
We use the same notation
$\rho$
for the induced metric on
$S^n\sub\R^{n+1}$,
and for the induced metric on
$\wh\R^n=\{x_{n+1}=0\}\cup\{\infty\}\sub\wh\R^{n+1}$.

We can generalize this classical situation
to the infinite dimensional case, by
replacing 
$\R^n$ by $\ell^2$, and
$\R^{n+1}$ by
$\R \times \ell^2$ with an additional $0$-coordinate.
Then the unit sphere $S^{\infty}$ in 
$\R \times \ell^2$ is the boundary of the infinite dimensional
hyperbolic space, where the Bourdon metric
(with respect to the origin) is the metric
$\frac{1}{2} \rho$.
The map $\pi: \wh\ell^2\to S^{\infty}$ restricted to the bounded
subset
$g(Z) \subset \ell^2$ is bi-Lipschitz.
Thus
$f:Z \to S^{\infty}$,
$f= \pi\circ g$, is a snowflake map.

\begin{tabbing}

Thomas Foertsch,\hskip11em\relax \= Viktor Schroeder,\\

Institut f\"ur Mathematik \>
Institut f\"ur Mathematik \\

Universit\"at Bonn \>
Universit\"at Z\"urich\\

Beringstr. 1 \>
Winterthurer Strasse 190 \\

D-53111 Bonn \>  
CH-8057 Z\"urich, Switzerland\\

{\tt foertsch@math.uni-bonn.de}\> {\tt vschroed@math.unizh.ch}\\

\end{tabbing}

\end{document}